\documentclass[twoside]{article}
\usepackage[usenames,dvipsnames]{color}
\usepackage{graphicx}
\usepackage{camnum}
\usepackage{harvard}

\usepackage{amsmath}
\usepackage{eufrak}
    \usepackage{ulem}
    \usepackage{chngcntr}
\usepackage{mathtools}
\usepackage{bm}

\newtheorem{remark}{Remark}
\counterwithin{table}{section}

\usepackage{tikz,pgfplots}
        \pgfplotsset{compat = 1.3}
        \pgfplotsset{minor grid style={dotted}} \pgfplotsset{major grid
        style={dashed}}

        \pgfplotsset{every x tick label/.append style={font=\footnotesize,
        yshift=0.25ex}}

        \pgfplotsset{every y tick label/.append
        style={font=\footnotesize, xshift=0.25ex}}

        \definecolor{colorclassyorange}{rgb}{0.95000,0.32500,0.09800}
        \definecolor{colorchromeyellow}{rgb}{1.00000,0.6549,0}%
        \definecolor{colorpaleyellow}{rgb}{1.00000,0.8549,0.1}%

        \definecolor{colorclassyblue}{rgb}{0.00000,0.44706,0.74118}%
        \definecolor{colorpurple}{rgb}{0.49400,0.18400,0.55600}%
        \definecolor{colorfuschia}{rgb}{0.95039,0.0,0.95039}%
        \definecolor{colorlemongreen}{rgb}{0.6,0.8,0}%
        \definecolor{colorreal}{rgb}{0.92941,0.79412,0.12549}%
        \definecolor{colorimag}{rgb}{0.00000,0.49804,0.00000}%
        \definecolor{colorabs}{rgb}{1.00000,0.00000,0.00000}%

\newcommand{\ee}{{\mathrm e}}
\newcommand{\ii}{{\mathrm i}}

\newcommand{\FFF}{\mathcal{F}}
\newcommand{\AAA}{\mathcal{A}}
\newcommand{\DDD}{\mathcal{D}}
\newcommand{\WW}{\mathcal{W}}
\newcommand{\dx}{\partial_x}
\newcommand{\dy}{\partial_y}

\def\O#1{{\cal O}\left(#1\right)}

\newcommand{\Int}[4]{\int_{#2}^{#3}\!#4\,\mathrm{d}#1}

\newcommand{\ve}{\varepsilon}
\newcommand{\schr}{Schr\"odinger }

\newcommand{\norm}[2]{\left\| #2 \right\|_{#1}}

\numberwithin{equation}{section}

\begin{document}

\title{Compact schemes for laser-matter interaction in \schr equation}

\author{Arieh Iserles,\footnote{Department of Applied Mathematics and Theoretical Physics,
 University of Cambridge, Wilberforce Rd, Cambridge CB3 0WA, UK.}\ \
  Karolina Kropielnicka\footnote{Institute of Mathematics, Polish Academy of Sciences,
  8 \'Sniadeckich Street, 00-656 Warsaw, Poland.}\ \  \&
  Pranav Singh\footnote{Mathematical Institute, Andrew Wiles Building, University of Oxford, Radcliffe Observatory Quarter, Woodstock Rd, Oxford OX2 6GG, UK and Trinity College, University of Oxford, Broad Street, Oxford OX1 3BH, UK.}}
\maketitle

\vspace*{6pt}

\vspace*{6pt} \noindent   \textbf{Keywords:} Schr\"odinger equation, laser
potential, higher order method, compact methods, splitting methods, time
dependent potentials,  Magnus expansion

\begin{abstract}
Numerical solutions for laser-matter interaction in \schr equation has many
applications in theoretical chemistry, quantum physics and condensed matter
physics. In this paper we introduce a methodology which allows, with a
small cost, to extend any fourth-order scheme for \schr equation with
time-indepedent potential to a fourth-order method
for \schr equation with laser potential. 
These fourth-order methods improve upon many leading schemes of order six
due to their low costs and small error constants.




%

\end{abstract}


\setcounter{section}{0}

\section{Introduction}
\label{sec:intro}

In this paper we are concerned with developing highly efficient numerical approaches for laser-matter interaction in the \schr equation,
\begin{equation}\label{eq:schr}
\ii \ve \partial_t u(\MM{x},t)=\big [-\ve^2 \Delta + V_0(\MM{x}) + \MM{e}(t)^\top  \MM{x}\big] u(\MM{x},t),\ u(\MM{x},0)=u_0(\MM{x}),
\end{equation}
where $\ t \geq 0$, $\MM{x}=(x_1,\ldots,x_n)\in\BB{R}^n$ and the laser term
$\MM{e}(t)=(e_1(t),\ldots,e_n(t))$ is an $\BB{R}^n$ valued function of $t$.
If the direction of the laser is assumed to be fixed, $\MM{e}(t)$ is of the
form
\[ \MM{e}(t) = e(t) \hat{\MM{\mu}},\quad e(t)\in \BB{R},\ \hat{\MM{\mu}} \in \BB{R}^n,\ \norm{2}{\hat{\MM{\mu}}} =1, \]
A very specific case, which is frequently used is
$\hat{\MM{\mu}}=(1,0,\ldots,0)$, so that
\[ V(\MM{x},t) = V_0(\MM{x}) - e(t) x_1.\]
The parameter $\ve$ in \R{eq:schr} acts like Planck's constant. This
parameter is $1$ when working in the atomic scaling and is very small, $0
<\ve \ll 1$, when working in the semiclassical regime. The methods developed
in this paper are equally effective in both regimes and will cover the
special cases of lasers as well.

\schr equations under the influence of lasers play a significant role in
quantum physics and theoretical chemistry, as they aid in the simulation and
design of systems and processes at atomic and molecular scales. Highly
accurate and cost effective schemes are required, for instance, in
applications such as optimal control strategies for shaping lasers where the
numerical solutions for these equations are used repeatedly within an
optimisation routine \cite{AmstrupChirped,MeyerOptimal,CoudertOptimal}.

The presence of time-dependant potentials in the \schr equation makes the
challenge of designing an efficient method significantly harder than the case
of time-independent potentials since typical strategies for this case require
utilisation of Magnus expansion at each time step, which involves nested
integrals of nested commutators of large matrices.

The design of effective methods for time-dependent potentials is, therefore,
a very challenging and active research area in theoretical chemistry, quantum
physics and numerical mathematics. A wide range of methodologies has been
developed to effectively handle this case
\cite{TalEzer1992,Peskin,SanzSerna96,TremblayCarrington04,KlaDimBrig2009,NdongEzerKosloff10,alvermann2011hocm,BlanesCasasMurua17,blanes17quasimagnus,Schaefer,IKS18sinum,IKS18jcp}.
The aforementioned methods can indeed be applied to (\ref{eq:schr}). However,
apart the third-order method from \cite{KlaDimBrig2009}, they are not
specialised for the case of lasers, where the structure of the potential can
be exploited to yield more efficient methods.

%

In this paper we propose a fourth-order numerical method, highly specialised
for the \schr equation under the influence of a laser (\ref{eq:schr}). The
main merits of the proposed method are its low costs and the ease with which
existing fourth-order implementations for the time-independent potentials can
be extended to handle time-dependent potentials. The cost of the proposed
methods is only marginally higher than the fastest fourth-order methods
dedicated for \schr equation with time-independent potentials, $V_0(\MM{x})$.

The two main ingredients in the derivation of our method are (i) the
simplification of commutators in the Magnus expansion by exploiting the
special form of laser potentials (ii) approximating the exponential of this
Magnus expansion by appropriate fourth-order splittings, again exploiting the
special form of the potential. In step (ii), we face a choice -- we may opt
for (a) a combination of Strang splitting with a fourth-order scheme for
time-independent potentials or (b) a compact splitting featuring positive
coefficients directly for this Magnus expansion.

In the application of Strang splitting in the first case, we extract the
smallest component of the Magnus expansion as the outer term (with a cost of
only two additional FFTs, due to its special structure). This yields fourth
order accuracy. The inner term can then be approximated with any of the many
suitable methods for \schr equation with time-independent potentials. This
may be the highly optimised splitting of \cite{blanesandmoan}, the compact
splitting of \cite{ChinChen02} (if the gradient of potential is available) or
any other fourth-order method, depending on the requirements or availability
of an existing implementation.


%
%


In Section~\ref{sec:general} we obtain, via simplification of commutators,
the optimal form of Magnus expansion, which can be handled efficiently. In
Section~\ref{sec:StrangOnMagnus} we present application of Strang splitting,
explaining why it serves here as a fourth order method. In its subsections,
in turn, we describe two (out of many) options that could be applied for the
exponentiation of the inner exponential. We also briefly discuss the
computational aspects and explain the low cost of this approach.
Section~\ref{sec:chin4mag} pursues the second alternative, presenting a
compact splitting scheme applied directly to the Magnus expansion of
Section~\ref{sec:general}. Numerical examples are described in
Section~\ref{sec:numerics}, while our conclusions are summarised in
Section~\ref{sec:conclusions}.

\section{Application of the Magnus expansion}
\label{sec:general}
The \schr equation \R{eq:schr} can be rewritten in the form
\begin{equation}\label{eq:schr2}
\partial_t u(\MM{x},t)= \AAA(\MM{x},t) u(\MM{x},t),\quad  \MM{x} \in  \BB{R}^n,\ t \geq 0,\ u(\MM{x},0)=u_0(\MM{x}),
\end{equation}
where $\AAA(\MM{x},t) = \ii \ve \Delta -\ii \ve^{-1} \left(
V_0(\MM{x})+\MM{e}(t)^\top \MM{x}\right)$. A fourth-order numerical scheme
for this equation can be obtained by resorting to a Magnus expansion
\cite{magnus54},
\[ \MM{u}(t+h) = \ee^{\Theta_2(t+h,t)} \MM{u}(t),\]
where the Magnus series has been truncated to the first two
terms\footnote{Commutator of the operators $\AAA$ and $\mathcal{B}$ is
defined here as $[\AAA,\mathcal{B}]=\AAA \circ \mathcal{B} - \mathcal{B}\circ
\AAA$.},
\[\Theta_2(t+h,t) = \Int{\zeta}{0}{h}{\AAA(t+\zeta)} - \Frac12
\Int{\zeta}{0}{h}{\Int{\xi}{0}{\zeta}{ [\AAA(t+\xi),A(t+\zeta)]  }   }.\]
\[\]

It is a simple matter to verify that
\begin{Eqnarray*}
\nonumber \Theta_2(t+h,t)&=& \ii h \ve \Delta - \ii h \ve^{-1} V_0(\MM{x}) - \ii
\ve^{-1}
\left(\Int{\zeta}{0}{h}{\MM{e}(t+\zeta)}\right)^\top  \MM{x}\\
\nonumber & & - \Frac12 \left(\Int{\zeta}{0}{h}{\Int{\xi}{0}{\zeta}{
\left[\MM{e}(t+\zeta) - \MM{e}(t+\xi) \right] } }\right)^\top [\Delta,
\MM{x}].
\end{Eqnarray*}
We note that
\[ [\Delta,\MM{x}]u = \sum_{j=1}^n \left( \partial_{x_j}^2 (\MM{x} u) - \MM{x} \partial_{x_j}^2 u \right) =  2 \nabla u, \]
and the fact that the integral over the triangle in $\Theta_2$ can be
rewritten as a univariate integral,
\[\Int{\zeta}{0}{h}{\Int{\xi}{0}{\zeta}{ \left[\MM{e}(t+\zeta) -
\MM{e}(t+\xi) \right] } } = 2 \Int{\zeta}{0}{h}{ \left(\zeta -
\Frac{h}{2}\right)\MM{e}(t+\zeta)}.\] Using these relations, we recast the Magnus expansion up to order four in the form
\begin{equation}
\label{eq:Mag} \Theta_2(t+h,t) = \ii h \ve \Delta - \ii h \ve^{-1} \widetilde{V}(\MM{x},t,h)  - \MM{s}(t,h)^\top  \nabla,
\end{equation}
where
\begin{align}
\label{eq:Vtilde}
\widetilde{V}(\MM{x},t,h) &= V_0(\MM{x}) + \MM{r}(t,h)^\top  \MM{x},\\
\label{eq:r}\MM{r}(t,h) &= \Frac{1}{h} \Int{\zeta}{0}{h}{\MM{e}(t+\zeta)},\\
\label{eq:s}\MM{s}(t,h) &= 2 \Int{\zeta}{0}{h}{
\left(\zeta - \Frac{h}{2}\right)\MM{e}(t+\zeta)}.
\end{align}
The numerical exponentiation of (\ref{eq:Mag}) is incredibly costly unless it
is split properly. In Sections~\ref{sec:splitmag} and \ref{sec:chin4mag}, we
will explore two fourth-order exponential splitting strategies for
approximating this exponential.


\section{A Strang splitting of the Magnus expansion}
\label{sec:StrangOnMagnus}
\label{sec:splitmag} A Strang splitting of $\Theta_2(t+h,t) = X + Y$, with $X
= \ii h \ve \Delta - \ii h \ve^{-1} \widetilde{V}(\MM{x},t,h)$ and $Y
=-\MM{s}(t,h)^\top  \nabla$,
\begin{equation}
\label{eq:strang}
\exp\left(-\Frac12 \MM{s}(t,h)^\top  \nabla\right) \ee^{\ii h \ve \Delta -
\ii  h\ve^{-1} \widetilde{V}(\MM{x},t,h) } \exp\left(-\Frac12 \MM{s}(t,h)^\top
 \nabla\right), \end{equation} is an order-four method, not order two as
usually expected  of  Strang splitting.

To see this, note that according
to the symmetric Baker--Campbell--Hausdorff (sBCH) formula,
\[ \ee^{\frac12 Y} \ee^{X} \ee^{\frac12 Y} =
\ee^{\mathrm{sBCH}(Y,X)}, \ \ \mathrm{sBCH(X, Y)} = X+Y -
\left(\Frac{1}{24}[[Y,X],X] + \Frac{1}{12} [[Y,X],Y] \right) +
\mathrm{h.o.t.}, \] the Strang splitting approximates the exponential of
$X+Y$ up to the commutators $\Frac{1}{24}[[Y,X],X] +
\Frac{1}{12} [[Y,X],Y]$. However, unlike the usual application of Strang splitting, where $X,Y=\O{h}$,
 the second term in the Magnus expansion, $Y = - \MM{s}(t,h)^\top  \nabla$,
scales as $\O{h^3}$ \cite{iserles00lgm}. Consequently, the largest commutator in our case is $\Frac{1}{24}[[Y,X],X] = \O{h^5}$, and the error in the
application of Strang splitting used for deriving \R{eq:strang} is $\O{h^5}$.

The central exponent in \R{eq:strang} can now be approximated via any
order-four method for time-independent potentials. In other words, as long as
an effective integrator for time-independent potentials is available, it need
only be modified in the following manner in order to convert it to an
order-four scheme for time-dependent potentials:

\begin{enumerate}
\item In each step $\widetilde{V}(\MM{x},t,h)$ has to be recomputed.
    However, as we can see from \R{eq:Vtilde}, this only requires the
    computation of the time-integral of the laser pulse \R{eq:r}, and not
    the re-evaluation of $V_0$, which needs to be computed only once.
\item The potential $\widetilde{V}$ should be used in place of $V_0$ in the
    existing fourth-order scheme for time-independent potentials.
\item Lastly, we need to compute the outermost exponentials,
    $\exp\left(-\Frac12 \MM{s}(t,h)^\top  \nabla\right)$. Note that this
    can be combined across two consecutive steps of \R{eq:strang},
    \[\exp\left(-\Frac12 \MM{s}(t+h,h)^\top
    \nabla\right) \exp\left(-\Frac12 \MM{s}(t,h)^\top  \nabla\right) =
    \exp\left(-\Frac12 \left[\MM{s}(t+h,h) + \MM{s}(t,h) \right]^\top
    \nabla\right),\] so that the overall scheme only involves one
    additional exponential per step.
\end{enumerate}
%

In general, the evaluation of the extra exponential
should be no more expensive than the exponential of the Laplacian,
which is routinely employed in schemes for the \schr equation with time-independent potentials.
In some cases it is possible to combine this exponential in such a
manner that no additional cost is incurred in the scheme, as we shall see in Section~\ref{sec:combine}.

In the following sections we will consider two concrete examples of
order-four methods for time-independent potentials for approximation of the
central exponential of \R{eq:strang}.

\subsection{Approximation of the inner term}

A wide range of exponential splittings are readily available and easily implementable
for the approximation of the \schr equation with time-independent potentials \cite{mclachlan02sm,blanes06sso,blanes08sac}.

The choice of appropriate exponential splittings may be governed by a need
for fewer exponentials, lower error constants or other constraints such as
positivity of coefficients. The examples we will consider are the highly
optimised order-four exponential splitting of \cite{blanesandmoan} and the
compact splitting scheme from \cite{ChinChen02}.

\subsubsection{Classical splittings}
\label{sec:blanes} The optimised splitting from \cite{blanesandmoan} for
$X,Y=\O{h}$,
\begin{equation}\label{eq:blanes}
  \ee^{X+Y} =  \ee^{a_1 X} \ee^{b_1 Y}  \ee^{a_2 X} \ee^{b_2 Y} \ee^{a_3 X} \ee^{b_3 Y}  \ee^{a_4 X}  \ee^{b_3 Y}  \ee^{a_3 X} \ee^{b_2 Y} \ee^{a_2 X} \ee^{b_1 Y}  \ee^{a_1 X} + \O{h^5},
\end{equation}
where
\begin{align*}
a_1 &= 0.0792036964311957, \quad & b_1 &= 0.209515106613362,\\
a_2 &= 0.353172906049774, \quad & b_2 &= -0.143851773179818,\\
a_3 &= -0.0420650803577195, \quad & b_3 &= 1/2 - b_1 - b_2,\\
a_4 &= 1 - 2(a_1 + a_2 +a_3), & &
\end{align*}
is known to have a very small error constant \cite{mclachlan02sm}. The
choices
\[ X = \ii h \ve \Delta = \O{h}, \quad Y = - \ii h \ve^{-1} \widetilde{V}(\MM{x},t,h) = \O{h},\]
result in an order-four splitting for the central exponent of \R{eq:strang}.
Combining this with the outer exponents in \R{eq:strang} completes one
example of a fourth-order scheme for time-dependent potentials.

\subsubsection{Compact splittings}
\label{sec:chin4}
Another concrete example results from using \cite{ChinChen02} for the approximation of the central exponential,
\begin{equation}\label{eq:chin4}
\ee^{X+Y} =
\ee^{\frac{1}{6}Y}\ee^{\frac{1}{2}X}\ee^{\frac{2}{3}(Y+\frac{1}{48}[[X,Y],Y])}\ee^{\frac{1}{2}X}\ee^{\frac{1}{6}Y}+\O{h^5}.
\end{equation}
In the case of
\[ X = \ii h \ve \Delta = \O{h}, \quad Y = - \ii h \ve^{-1} \widetilde{V}(\MM{x},t,h) = \O{h},\]
while the nested commutator reduces to a function,
\[ [[X,Y],Y] = -\ii h^3 \ve^{-1} [[\Delta,  \widetilde{V} ], \widetilde{V}] = -2 \ii h^3 \ve^{-1} \sum_{j=0}^n [(\partial_{x_j} \widetilde{V}) \partial_{x_j},\widetilde{V}] = -2\ii h^3 \ve^{-1} \nabla \widetilde{V}^\top \nabla \widetilde{V}. \]
Thus, \R{eq:chin4} reduces to
\begin{equation}\label{eq:chin4concrete} \ee^{ -\frac{1}{6} \ii h \ve^{-1} \widetilde{V}} \ee^{\frac{1}{2} \ii h \ve \Delta} \ee^{-\frac{2}{3}  \ii h \widehat{V}} \ee^{\frac{1}{2} \ii h \ve \Delta}\ee^{-\frac{1}{6} \ii h \ve^{-1} \widetilde{V}},
\end{equation}
where $\widehat{V}$ is an $\O{h^3}$ perturbation of $\widetilde{V}$,
\begin{equation} \label{eq:Vc} \widehat{V} = \widetilde{V} -\Frac{h^2}{24} (\nabla \widetilde{V}) ^\top (\nabla \widetilde{V}).
\end{equation}
The overall order-four scheme is obtained by substituting the central
exponent in \R{eq:strang} with \R{eq:chin4concrete},
\begin{equation}
\label{eq:MaStCC}
\ee^{-\Frac12 \MM{s}(t,h)^\top \nabla} \ee^{ -\frac{1}{6} \ii h \ve^{-1}
\widetilde{V}} \ee^{\frac{1}{2} \ii h \ve \Delta} \ee^{-\frac{2}{3}  \ii h
\ve^{-1} \widehat{V}} \ee^{\frac{1}{2} \ii h \ve \Delta}\ee^{-\frac{1}{6} \ii
h \ve^{-1} \widetilde{V}} \ee^{-\Frac12 \MM{s}(t,h)^\top \nabla}.
\end{equation}
This scheme involves the evaluation of $\nabla \widetilde{V}$. However, the typically expensive part, $\nabla V_0$, needs to be computed only
once since
\[ \nabla \widetilde{V} = \nabla V_0 + \MM{r}(t,h). \]

\subsubsection{Combining exponentials}
\label{sec:combine}
A further optimisation is possible in some cases. Upon replacing the
inner exponential in \R{eq:strang} by \R{eq:blanes} in Section~\ref{sec:blanes}, it should be observed that
the outermost exponentials commute and, therefore, can be combined,
\[ \exp\left(-\Frac12 \MM{s}(t,h)^\top
    \nabla\right) \ee^{ \ii a_1 h \ve \Delta}= \exp\left(\ii a_1 h \ve \Delta -\Frac12 \MM{s}(t,h)^\top
    \nabla \right).\]
In practice, this combined exponential is often not much harder to compute than the exponential of the Laplacian.

For instance, under spectral collocation the differentiation matrices are circulant and are diagonalised via
Fourier transforms\footnote{We use $\leadsto$ to denote discretisation.},
\[ \partial_{x}^k \leadsto \MM{D}_{k,x} = \FFF_{x}^{-1} \DDD_{c_{k,x}} \FFF_{x}, \]
where $\DDD_{c_{k,x}}$ is a diagonal matrix and the values along its diagonal, $c_{k,x}$, comprise the symbol of the
$k$th differentiation matrix, $\MM{D}_{k,x}$. In two dimensions the exponential of the Laplacian term, $\ii a_1 h \ve \Delta$, alone is implemented as
\[ \ee^{\ii a_1 h \ve \Delta} v \leadsto \FFF_{x}^{-1} \DDD_{\exp(\ii a_1 h \ve c_{2,x})} \FFF_{x} \FFF_{y}^{-1} \DDD_{\exp(\ii a_1 h \ve c_{2,y})} \FFF_{y} \MM{v}, \]
using four Fast Fourier Transforms (FFTs).
Similarly,
\[\ee^{\ii a_1 h \ve \Delta -\frac12 \MM{s}(t,h)^\top
    \nabla} =
\ee^{\ii h \ve a_1 \partial_{x}^2 - \frac12 s_x \partial_{x}}  \ee^{\ii h \ve a_1 \partial_{y}^2 - \frac12 s_y \partial_{y}}, \]
where $\MM{s} = (s_x,s_y)$, can be implemented as
\begin{eqnarray*}
&&\FFF_x^{-1} \DDD_{\exp( \ii h \ve a_1 c_{2,x} - \frac12
s_x c_{1,x})}
 \FFF_x \FFF_y^{-1} \DDD_{\exp( \ii h \ve a_1 c_{2,x} - \frac12
s_y c_{1,x})}
 \FFF_y,
\end{eqnarray*}
using the same number of FFTs\footnote{Here $c_{1,x},c_1{1,y},c_{2,x}$ and $c_{2,y}$ are the symbols of the differentiation matrices
corresponding to $\dx,\dy,\dx^2$ and $\dy^2$, respectively.}.

In this way, the proposed method for time-dependent potentials of the form
$V_0 + \MM{e}(t)^\top \MM{x}$ carries no additional expense compared to the
order-four splitting of \cite{blanesandmoan} for $V_0$ alone.

\begin{remark}
Note that another variant of the splitting can be obtained by swapping the
choices of $A$ and $B$. In this variant the outermost exponential of
\R{eq:blanes} is not a Laplacian and, therefore, cannot be combined in this
way with the outermost exponential of \R{eq:strang}. In this case the
proposed method requires one additional exponential.
\end{remark}

\begin{remark}
An alternative, of course, is to perform any order-four splitting on the
Magnus expansion {$\Theta_2(t+h,t)$} directly by choosing $A = \ii h \ve
\Delta - \frac12 \MM{s}^\top \nabla$ and $B=- \ii h \ve^{-1}
\widetilde{V}(\MM{x},t,h)$. In light of the comments in this section, such
variants also do not require additional expense in the case of discretisation
via spectral collocation, for instance.
\end{remark}


\section{A compact splitting of the Magnus expansion}
\label{sec:chin4mag}
As we have noted previously in subsecton~\ref{sec:chin4}, the highly efficient compact splitting scheme of \cite{ChinChen02},
\begin{equation}\tag{\ref{eq:chin4}}
\ee^{X+Y} =
\ee^{\frac{1}{6}Y}\ee^{\frac{1}{2}X}\ee^{\frac{2}{3}(Y+\frac{1}{48}[[X,Y],Y])}\ee^{\frac{1}{2}X}\ee^{\frac{1}{6}Y}+\O{h^5}.
\end{equation}
for $X,Y=\O{h}$, features a central exponent which, in principle, has a nested
commutator of $X$ and $Y$. In practice, however, when $X$ is the Laplacian term, this nested commutator reduces to a function, as we have seen in subsection~\ref{sec:chin4}. This results in a central
exponent which is an $\O{h^3}$ pertubation of the potential function.

Since the relation \R{eq:chin4} holds for arbitrary $X,Y$, in particular it also holds if
we use it for splitting $\Theta_2(t+h,t)$ into the two parts
\[X=\ii h \ve
\Delta -\MM{s}(t,h)^\top \nabla, \quad \mathrm{and}\quad Y=- \ii h \ve^{-1}
\widetilde{V}(\MM{x},t,h).\] Crucially, since the commutator of $\MM{s}^\top
\nabla$ with a function $f(x)$ reduces to another function\footnote{By noting
$[\dx, f] u = \dx (f u) - f (\dx u) = f \dx u + (\dx f) u - f (\dx u) = (\dx
f) u$, we conclude the operatorial identity $[\dx, f] = (\dx f)$.},
\[ [\MM{s}^\top \nabla,f] = \sum_{j=1}^n s_j [\partial_{x_j},f]
= \sum_{j=1}^n s_j \left(f \partial_{x_j} + (\partial_{x_j} f) - f
\partial_{x_j}\right) = \MM{s}^\top (\nabla f),\]
the nested commutator
\[ [[\MM{s}^\top \nabla,f],f] = [\MM{s}^\top (\nabla f),f] = 0,\]
vanishes.
 The central exponent in this splitting is, therefore, identical to the
one encountered in subsecton~\ref{sec:chin4},
\[ \Frac{2}{3}(Y+\Frac{1}{48}[[X,Y],Y]) = -\Frac{2}{3} \ii h \ve^{-1} \widehat{V}
= -\Frac{2}{3} \ii h \ve^{-1} (\widetilde{V}-\Frac{h^2}{24} (\nabla \widetilde{V})^\top (\nabla \widetilde{V})).\]

Thus, a fourth-order compact splitting,
\begin{equation} \label{eq:MaCC}
\ee^{ -\frac{1}{6} \ii h \ve^{-1} \widetilde{V}} \ee^{\frac{1}{2} \ii h \ve
\Delta -\frac12 \MM{s}(t,h)^\top \nabla} \ee^{-\frac{2}{3}  \ii h \ve^{-1}
\widehat{V}} \ee^{\frac{1}{2} \ii h \ve \Delta-\frac12 \MM{s}(t,h)^\top
\nabla}\ee^{-\frac{1}{6} \ii h \ve^{-1} \widetilde{V}}, \end{equation} can be
directly implemented using five stages (exponentials) for the case of
time-dependent potentials.

This differs from the \cite{ChinChen02} for time-independent potentials only
in that a modified potential $\widetilde{V}$ is used in place of the
time-independent potential $V_0$ and a perturbation of the Laplacian term,
$\ii h \ve \Delta -\MM{s}(t,h)^\top \nabla,$ is used instead of the usual
$\ii h \ve \Delta$ occurring in \cite{ChinChen02}. In light of the remarks
made in Subsection~\ref{sec:combine}, neither of these add substantially to
the cost of the method.

%
%
%
%
%

\section{Numerical examples}
\label{sec:numerics}
In this section we will consider two one-dimensional numerical examples --
the first in atomic scaling, $\ve_1=1$, and the second in the semiclassical
regime of $\ve_2 = 10^{-2}$.

The initial conditions $u_{0,1}$ and $u_{0,2}$ for our numerical experiments
are Gaussian wavepackets (with zero initial momentum),
\[u_{0,k}(x) = (\delta_k  \pi)^{-1/4}  \exp\left( (-(x - x_0)^2 )/ (2 \delta_k)\right),\quad x_0 = -2.5,\ k=1,2,\]
with $\delta_1 = 0.2$ and $\delta_2 = 10^{-2}$ in the respective cases. These
wavepackets are sitting in the left well of the double well potentials,
\[V_{\mathrm{D}1}(x) = x^4 - 15 x^2 \quad \mathrm{and}\quad V_{\mathrm{D}2}(x) = \Frac15 x^4 - 2 x^2, \]
respectively, which act as the choice of $V_0$ in the two examples.

Our spatial domain is $[-10,10]$ and $[-5,5]$ in the two examples,
respectively, while the temporal domain is $[0,4]$ and $[0,\frac52]$,
respectively. We impose periodic boundaries on the spatial domains and resort
to spectral collocation for discretisation.

\begin{figure}[tbh]
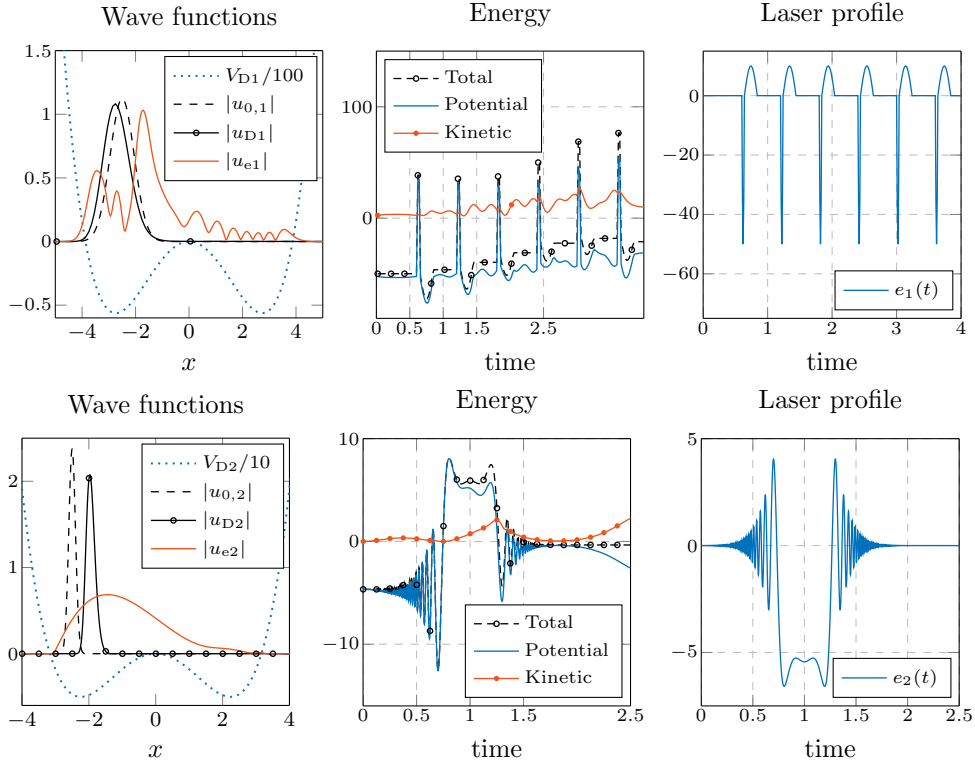

        \input{InitialFinalEx2.tex}
            \input{EnergyEx2.tex}
%
\definecolor{mycolor1}{rgb}{0.00000,0.44706,0.74118}%
\definecolor{mycolor2}{rgb}{0.95000,0.32500,0.09800}%
\definecolor{mycolor3}{rgb}{1.00000,0.99000,0.97000}%
\begin{tikzpicture}
\begin{axis}[%
width=1.35in, height=1.4in, at={(0in,0in)}, scale only axis, separate axis
lines,
xmin=0, xmax=4,
xlabel={time},
ymin=-75, ymax=15, title={Laser profile}, yticklabel style =
{font=\scriptsize,xshift=0.25ex},
xmajorgrids, ymajorgrids, xticklabel style =
{font=\scriptsize,yshift=0.25ex}, legend style={legend cell align=left,
align=left, draw=black,font=\scriptsize}, legend pos=south east] \addplot
[color=mycolor1, line width=0.5pt]
  table[row sep=crcr]{%
0	0\\
0.01	0\\
0.02	0\\
0.03	0\\
0.04	0\\
0.05	0\\
0.06	0\\
0.07	0\\
0.08	0\\
0.09	0\\
0.1	0\\
0.11	0\\
0.12	0\\
0.13	0\\
0.14	0\\
0.15	0\\
0.16	0\\
0.17	0\\
0.18	0\\
0.19	0\\
0.2	0\\
0.21	0\\
0.22	0\\
0.23	0\\
0.24	0\\
0.25	0\\
0.26	0\\
0.27	0\\
0.28	0\\
0.29	0\\
0.3	0\\
0.31	0\\
0.32	0\\
0.33	0\\
0.34	0\\
0.35	0\\
0.36	0\\
0.37	0\\
0.38	0\\
0.39	0\\
0.4	0\\
0.41	0\\
0.42	0\\
0.43	0\\
0.44	0\\
0.45	0\\
0.46	0\\
0.47	0\\
0.48	0\\
0.49	0\\
0.5	0\\
0.51	0\\
0.52	0\\
0.53	0\\
0.54	0\\
0.55	0\\
0.56	0\\
0.57	0\\
0.58	0\\
0.59	0\\
0.6	0\\
0.61	-35.3553390593273\\
0.62	-50\\
0.63	-35.3553390593274\\
0.64	0\\
0.65	1.56434465040231\\
0.66	3.09016994374948\\
0.67	4.53990499739547\\
0.68	5.87785252292474\\
0.69	7.07106781186548\\
0.7	8.09016994374948\\
0.71	8.91006524188368\\
0.72	9.51056516295153\\
0.73	9.87688340595138\\
0.74	10\\
0.75	9.87688340595138\\
0.76	9.51056516295154\\
0.77	8.91006524188368\\
0.78	8.09016994374948\\
0.79	7.07106781186547\\
0.8	5.87785252292473\\
0.81	4.53990499739546\\
0.82	3.09016994374947\\
0.83	1.5643446504023\\
0.84	1.01064309961486e-14\\
0.85	0\\
0.86	0\\
0.87	0\\
0.88	0\\
0.89	0\\
0.9	0\\
0.91	0\\
0.92	0\\
0.93	0\\
0.94	0\\
0.95	0\\
0.96	0\\
0.97	0\\
0.98	0\\
0.99	0\\
1	0\\
1.01	0\\
1.02	0\\
1.03	0\\
1.04	0\\
1.05	0\\
1.06	0\\
1.07	0\\
1.08	0\\
1.09	0\\
1.1	0\\
1.11	0\\
1.12	0\\
1.13	0\\
1.14	0\\
1.15	0\\
1.16	0\\
1.17	0\\
1.18	0\\
1.19	0\\
1.2	0\\
1.21	-35.3553390593273\\
1.22	-50\\
1.23	-35.3553390593274\\
1.24	0\\
1.25	1.56434465040231\\
1.26	3.09016994374948\\
1.27	4.53990499739547\\
1.28	5.87785252292474\\
1.29	7.07106781186548\\
1.3	8.09016994374948\\
1.31	8.91006524188368\\
1.32	9.51056516295154\\
1.33	9.87688340595138\\
1.34	10\\
1.35	9.87688340595138\\
1.36	9.51056516295153\\
1.37	8.91006524188367\\
1.38	8.09016994374946\\
1.39	7.07106781186546\\
1.4	5.87785252292471\\
1.41	4.53990499739548\\
1.42	3.09016994374948\\
1.43	1.56434465040232\\
1.44	1.01064309961486e-14\\
1.45	0\\
1.46	0\\
1.47	0\\
1.48	0\\
1.49	0\\
1.5	0\\
1.51	0\\
1.52	0\\
1.53	0\\
1.54	0\\
1.55	0\\
1.56	0\\
1.57	0\\
1.58	0\\
1.59	0\\
1.6	0\\
1.61	0\\
1.62	0\\
1.63	0\\
1.64	0\\
1.65	0\\
1.66	0\\
1.67	0\\
1.68	0\\
1.69	0\\
1.7	0\\
1.71	0\\
1.72	0\\
1.73	0\\
1.74	0\\
1.75	0\\
1.76	0\\
1.77	0\\
1.78	0\\
1.79	0\\
1.8	-7.3887043024834e-13\\
1.81	-35.3553390593273\\
1.82	-50\\
1.83	-35.3553390593274\\
1.84	3.48786849800863e-14\\
1.85	1.56434465040231\\
1.86	3.09016994374951\\
1.87	4.53990499739547\\
1.88	5.87785252292476\\
1.89	7.07106781186548\\
1.9	8.0901699437495\\
1.91	8.91006524188368\\
1.92	9.51056516295153\\
1.93	9.87688340595138\\
1.94	10\\
1.95	9.87688340595138\\
1.96	9.51056516295154\\
1.97	8.91006524188367\\
1.98	8.09016994374948\\
1.99	7.07106781186546\\
2	5.87785252292474\\
2.01	4.53990499739551\\
2.02	3.09016994374948\\
2.03	1.56434465040228\\
2.04	1.01064309961486e-14\\
2.05	0\\
2.06	0\\
2.07	0\\
2.08	0\\
2.09	0\\
2.1	0\\
2.11	0\\
2.12	0\\
2.13	0\\
2.14	0\\
2.15	0\\
2.16	0\\
2.17	0\\
2.18	0\\
2.19	0\\
2.2	0\\
2.21	0\\
2.22	0\\
2.23	0\\
2.24	0\\
2.25	0\\
2.26	0\\
2.27	0\\
2.28	0\\
2.29	0\\
2.3	0\\
2.31	0\\
2.32	0\\
2.33	0\\
2.34	0\\
2.35	0\\
2.36	0\\
2.37	0\\
2.38	0\\
2.39	0\\
2.4	0\\
2.41	-35.3553390593279\\
2.42	-50\\
2.43	-35.355339059328\\
2.44	0\\
2.45	1.56434465040234\\
2.46	3.09016994374948\\
2.47	4.53990499739544\\
2.48	5.87785252292474\\
2.49	7.0710678118655\\
2.5	8.09016994374948\\
2.51	8.91006524188367\\
2.52	9.51056516295154\\
2.53	9.87688340595139\\
2.54	10\\
2.55	9.87688340595138\\
2.56	9.51056516295153\\
2.57	8.91006524188366\\
2.58	8.09016994374946\\
2.59	7.07106781186549\\
2.6	5.87785252292477\\
2.61	4.53990499739548\\
2.62	3.09016994374945\\
2.63	1.56434465040232\\
2.64	4.1192675685653e-14\\
2.65	0\\
2.66	0\\
2.67	0\\
2.68	0\\
2.69	0\\
2.7	0\\
2.71	0\\
2.72	0\\
2.73	0\\
2.74	0\\
2.75	0\\
2.76	0\\
2.77	0\\
2.78	0\\
2.79	0\\
2.8	0\\
2.81	0\\
2.82	0\\
2.83	0\\
2.84	0\\
2.85	0\\
2.86	0\\
2.87	0\\
2.88	0\\
2.89	0\\
2.9	0\\
2.91	0\\
2.92	0\\
2.93	0\\
2.94	0\\
2.95	0\\
2.96	0\\
2.97	0\\
2.98	0\\
2.99	0\\
3	0\\
3.01	-35.3553390593267\\
3.02	-50\\
3.03	-35.3553390593268\\
3.04	0\\
3.05	1.56434465040228\\
3.06	3.09016994374948\\
3.07	4.53990499739544\\
3.08	5.87785252292474\\
3.09	7.07106781186546\\
3.1	8.09016994374948\\
3.11	8.91006524188367\\
3.12	9.51056516295154\\
3.13	9.87688340595137\\
3.14	10\\
3.15	9.87688340595138\\
3.16	9.51056516295153\\
3.17	8.91006524188369\\
3.18	8.0901699437495\\
3.19	7.07106781186549\\
3.2	5.87785252292471\\
3.21	4.53990499739548\\
3.22	3.09016994374952\\
3.23	1.56434465040232\\
3.24	0\\
3.25	0\\
3.26	0\\
3.27	0\\
3.28	0\\
3.29	0\\
3.3	0\\
3.31	0\\
3.32	0\\
3.33	0\\
3.34	0\\
3.35	0\\
3.36	0\\
3.37	0\\
3.38	0\\
3.39	0\\
3.4	0\\
3.41	0\\
3.42	0\\
3.43	0\\
3.44	0\\
3.45	0\\
3.46	0\\
3.47	0\\
3.48	0\\
3.49	0\\
3.5	0\\
3.51	0\\
3.52	0\\
3.53	0\\
3.54	0\\
3.55	0\\
3.56	0\\
3.57	0\\
3.58	0\\
3.59	0\\
3.6	0\\
3.61	-35.3553390593267\\
3.62	-50\\
3.63	-35.355339059328\\
3.64	0\\
3.65	1.56434465040228\\
3.66	3.09016994374948\\
3.67	4.53990499739544\\
3.68	5.87785252292474\\
3.69	7.07106781186546\\
3.7	8.09016994374948\\
3.71	8.91006524188367\\
3.72	9.51056516295152\\
3.73	9.87688340595137\\
3.74	10\\
3.75	9.87688340595138\\
3.76	9.51056516295155\\
3.77	8.91006524188369\\
3.78	8.0901699437495\\
3.79	7.07106781186549\\
3.8	5.87785252292477\\
3.81	4.53990499739548\\
3.82	3.09016994374952\\
3.83	1.56434465040232\\
3.84	4.1192675685653e-14\\
3.85	0\\
3.86	0\\
3.87	0\\
3.88	0\\
3.89	0\\
3.9	0\\
3.91	0\\
3.92	0\\
3.93	0\\
3.94	0\\
3.95	0\\
3.96	0\\
3.97	0\\
3.98	0\\
3.99	0\\
4	0\\
}; \addlegendentry{$e_1(t)$}

\end{axis}
\end{tikzpicture}%
\\
            \input{InitialFinalEx1.tex}
            \input{EnergyEx1.tex}
            \input{LaserProfileEx1.tex}
        \caption{[Example 1 (top row), Example 2 (bottom row)]. The initial condition
        $u_{0,k}$ evolves to $u_{\mathrm{D}k}$ at time $T$ under the influence of $V_0=V_{\mathrm{D}k}$ and
        to $u_{\mathrm{e}k}$ under $V_{\mathrm{e}k}=V_{\mathrm{D}k}+e_k(t) x$ (left); Evolution of energy under
        $V_{\mathrm{e}k}$ (centre); Laser profile $e_k(t)$ (right).
        Here, the potentials $V_{\mathrm{D}k}$ are scaled down for ease of presentation, $k=1$ for Example 1
        and $k=2$ for Example 2.}
        \label{fig:uIF}
    \end{figure}

When we allow the wave functions $u_{0,k}$ to evolve to $u_{\mathrm{D}k}$ at
the final time, $T$, under the influence of $V_{\mathrm{D}k}$ alone, they
remain largely confined to the left well (Figure~\ref{fig:uIF}, left column).

%
%
{\bf Time-dependent potential.} Superimposing a time dependent excitation of
the form $e(t)\,x$ on the potential -- used for modeling laser interaction --
we are able to exert control on the wave function. The time profile of the
laser used here is
\[ e_1(t) = \begin{cases} \sin(25 \pi t) & \quad t \in [\Frac35 n,\Frac35 n + \Frac{1}{25}],\quad n \geq 1,\\
                          \sin(5 \pi t) & \quad t \in (\Frac35 n + \Frac{1}{25},\Frac35 n + \Frac{6}{25}],\quad n \geq 1,
            \end{cases} \]
and
\[ e_2(t) = 10 \exp(-10 (t - 1)^2) \sin((500 (t - 1)^4 + 10)), \]
respectively. The former is a sequence of {\em asymmetric sine lobes} while
the latter is a highly oscillatory {\em chirped} pulse (Figure~\ref{fig:uIF},
right column). Such laser profiles are used routinely in laser control
\cite{AmstrupChirped}. Even more oscillatory electric fields often result
from optimal control algorithms \cite{MeyerOptimal,CoudertOptimal}.

The effective time-dependent potentials in the two examples are
\[ V_{\mathrm{e}k}(x,t) = V_{\mathrm{D}k}(x) + e_k(t) x, \quad k=1,2.\]
Under the influence of this time-dependent potential $V_{\mathrm{e}k}(x,t)$,
the initial wave-functions $u_{0,k}$ evolve to $u_{\mathrm{e}k}$, which are
not confined to the left well (Figure~\ref{fig:uIF}, left column). In this
case the total energy is not conserved (Figure~\ref{fig:uIF}, middle column).



{\bf Methods.} We denote the combination of the fourth-order splitting
\R{eq:blanes} with \R{eq:strang}, which is derived in
Subsection~\ref{sec:blanes}, by MaStBM (Magnus--Strang--Blanes--Moan). This
method implements the outermost exponential in \R{eq:strang} separately and
is expected to be the most costly of our schemes. As noted in
Subsection~\ref{sec:combine}, the outermost exponentials can be combined,
resulting in a method which will be labeled as MaStBMc (c for combined).

Another concrete example of our splittings is outlined in
Subsection~\ref{sec:chin4}, where a more efficient method \R{eq:MaStCC}
results by resorting to \cite{ChinChen02} for the approximation of the
central exponent in \R{eq:strang}. This scheme is denoted by MaStCC
(Magnus--Strang--Chin--Chen). It requires computation of an additional
exponential in comparison with the method of \cite{ChinChen02}. The most
efficient method in this class \R{eq:MaCC}, results from the direct
application of \cite{ChinChen02} to $\Theta_2$ in Section~\ref{sec:chin4mag},
and is denoted by MaCC (Magnus--Chin--Chen). This method should cost no more
than an implementation of the method from \cite{ChinChen02}, which, of
course, is  designed for time-independent potentials. In this way, the fourth
order schemes presented here are able to handle time-dependent potentials
with little or no additional costs.

To demonstrate the effectiveness of our fourth-order schemes, we will compare
them to a couple of sixth-order methods. The first of these is the
sixth-order optimised method CF6:5Opt proposed in \cite{alvermann2011hocm},
which is labeled as AF in this section. AF is accompanied by a postfix which
refers to the number of Lanczos iterations (5, 10 or 50) used in the method.
The second scheme is the sixth-order method of \cite{IKS18jcp}, denoted by
IKS, which is specialised for the semiclassical regime. In all cases
presented here IKS, utilises three Lanczos iterations for the exponentiation
of its first non-trivial exponent, $W^{[2]}$, and two Lanczos iterations for
the exponentiation of the innermost exponent, $\WW^{[3]}$.

The first of these methods, AF, is not specialised  for laser potentials,
highly oscillatory potentials or for the semiclassical regime. The second,
IKS, is optimised for highly oscillatory potentials and does provide certain
optimisations for the case of lasers, which have been employed in these
experiments. However, it is not designed for atomic scaling.

In the first example, the integrals in our schemes, \R{eq:r} and \R{eq:s}, as
well as those in IKS are discretised via three Gauss--Legendre knots, while
in the second case eleven Gauss--Legendre knots are used in order to
adequately resolve the highly oscillatory potential. Not discretising these
integrals at the outset allows us flexibility in deciding a quadrature
strategy at the very end, as discussed by \cite{IKS18sinum,IKS18jcp}.
Effectively, this strategy allows us to resolve a highly oscillatory
potential despite using large time steps in the propagation of the solution.
In contrast, the method AF uses a fixed number (four) of Gauss--Legendre
knots for all cases.

\begin{figure}[tbh]
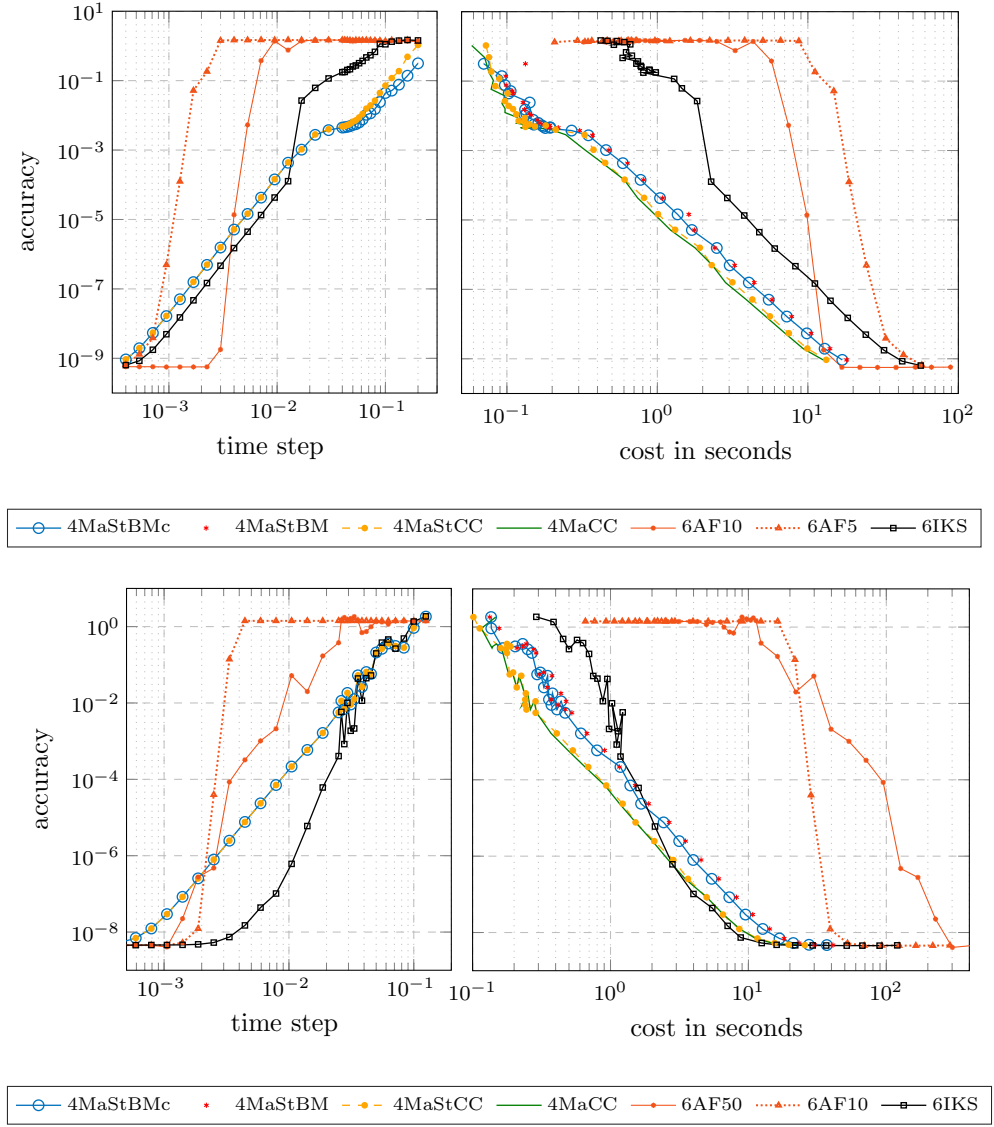

    \centering
%
%
\definecolor{mycolor1}{rgb}{0.00000,0.44700,0.74100}%
\definecolor{mycolor2}{rgb}{0.85000,0.32500,0.09800}%
\definecolor{mycolor3}{rgb}{0.92900,0.69400,0.12500}%
\definecolor{mycolor4}{rgb}{0.49400,0.18400,0.55600}%
\definecolor{mycolor5}{rgb}{0.46600,0.67400,0.18800}%
\definecolor{mycolor6}{rgb}{0.30100,0.74500,0.93300}%
\definecolor{mycolor7}{rgb}{0.63500,0.07800,0.18400}%
\begin{tikzpicture}

\begin{axis}[%
width=1.7in, height=2in, scale only axis, xmode=log, xmin=0.0003, xmax=0.3,
xminorticks=true, ymode=log, ymin=1e-10, ymax=10, yminorticks=true, axis
background/.style={fill=white}, xticklabel style =
{xshift=-0.25ex,yshift=-0.1ex}, ylabel={accuracy}, xlabel={time step
\mbox{}}, grid=both]

\addplot [color=colorclassyblue, line width=0.5pt,solid, mark=o, mark
options={solid, colorclassyblue}]
  table[row sep=crcr]{%
0.2	0.315073971584099\\
0.16	0.138656522753413\\
0.133333333333333	0.0766291104390896\\
0.114285714285714	0.0515648726624736\\
0.1	0.0441490382603395\\
0.0888888888888889	0.0241499327978784\\
0.08	0.0151479401677304\\
0.0727272727272727	0.011452266307402\\
0.0666666666666667	0.0103988241843331\\
0.0615384615384615	0.00763294863939209\\
0.0571428571428571	0.0062437019612876\\
0.0533333333333333	0.00556255171418981\\
0.05	0.00529064217606149\\
0.0470588235294118	0.00487089338911341\\
0.0444444444444444	0.00466497250515616\\
0.0421052631578947	0.00453274168827438\\
0.04	0.00444878125418753\\
0.0298507462686567	0.0037531573445756\\
0.0224719101123596	0.00274205401320506\\
0.0168067226890756	0.00102747561947899\\
0.0126182965299685	0.000431689965920295\\
0.00947867298578199	0.000142746190342925\\
0.00710479573712256	4.25722872739458e-05\\
0.00533333333333333	1.44530322492006e-05\\
0.004	5.12667649628139e-06\\
0.00299850074962519	1.55593059813843e-06\\
0.00224845418774592	4.91517706132474e-07\\
0.00168634064080944	1.57108617874099e-07\\
0.00126462219411951	5.03442351814988e-08\\
0.000948541617263457	1.64505765628687e-08\\
0.000711237553342817	5.39685056841894e-09\\
0.00053340445392719	1.95198135878671e-09\\
0.0004	9.28393959914495e-10\\
};

\addplot [color=colorchromeyellow,line width=0.5pt, dashed, mark=*,mark
size=1,mark options={solid},mark repeat=1]
  table[row sep=crcr]{%
0.2	1.061617989908\\
0.16	0.486154601499711\\
0.133333333333333	0.188084428891699\\
0.114285714285714	0.118230161358153\\
0.1	0.0720423805257222\\
0.0888888888888889	0.0442786668479905\\
0.08	0.0266119405061286\\
0.0727272727272727	0.0193414019339609\\
0.0666666666666667	0.0156454518101673\\
0.0615384615384615	0.0110289124533264\\
0.0571428571428571	0.0086854733959855\\
0.0533333333333333	0.00716199808259354\\
0.05	0.00658130004670823\\
0.0470588235294118	0.00567306756885257\\
0.0444444444444444	0.00525948241169924\\
0.0421052631578947	0.00492511050876687\\
0.04	0.00478400316871154\\
0.0298507462686567	0.00398833165920196\\
0.0224719101123596	0.00281789351733468\\
0.0168067226890756	0.00105118398308029\\
0.0126182965299685	0.000439266950282615\\
0.00947867298578199	0.000145176358134669\\
0.00710479573712256	4.33302670178804e-05\\
0.00533333333333333	1.46946975666681e-05\\
0.004	5.20360478905995e-06\\
0.00299850074962519	1.58019381740511e-06\\
0.00224845418774592	4.99185096630157e-07\\
0.00168634064080944	1.59530778953042e-07\\
0.00126462219411951	5.11093193470278e-08\\
0.000948541617263457	1.66915179551255e-08\\
0.000711237553342817	5.47159235498498e-09\\
0.00053340445392719	1.97391563166137e-09\\
0.0004	9.33793917905122e-10\\
};

\addplot [color=colorclassyorange, solid, mark=10-pointed star,mark
size=1,mark options={solid},mark repeat=1]
  table[row sep=crcr]{%
0.2	1.33821393223762\\
0.16	1.37323113229344\\
0.133333333333333	1.3331821823505\\
0.114285714285714	1.436849894784\\
0.1	1.39070636715556\\
0.0888888888888889	1.42183042187019\\
0.08	1.44124920796303\\
0.0727272727272727	1.45861287413695\\
0.0666666666666667	1.45267146485144\\
0.0615384615384615	1.46812397669847\\
0.0571428571428571	1.47139924554659\\
0.0533333333333333	1.43721064143572\\
0.05	1.43093588612415\\
0.0470588235294118	1.44077736942454\\
0.0444444444444444	1.44288357096747\\
0.0421052631578947	1.46219491157887\\
0.04	1.45660029318037\\
0.0298507462686567	1.48935047107182\\
0.0224719101123596	1.48739834666741\\
0.0168067226890756	1.38509745763673\\
0.0126182965299685	0.757774469514455\\
0.00947867298578199	1.38028293450082\\
0.00710479573712256	0.377652912207918\\
0.00533333333333333	0.00530678441546813\\
0.004	1.36820334923103e-05\\
0.00299850074962519	1.79446508659599e-09\\
0.00224845418774592	5.68608893649509e-10\\
0.00168634064080944	5.65192067566052e-10\\
0.00126462219411951	5.63776255412919e-10\\
0.000948541617263457	5.62985289323943e-10\\
0.000711237553342817	5.72402926201176e-10\\
0.00053340445392719	5.72676812235436e-10\\
0.0004	5.74625424583439e-10\\
};

\addplot [color=colorclassyorange, line width=0.75pt, densely dotted,
mark=triangle,mark size=1,mark options={solid},mark repeat=1]
  table[row sep=crcr]{%
0.2	1.32539816793023\\
0.16	1.46859467707918\\
0.133333333333333	1.36351262541305\\
0.114285714285714	1.45776970978599\\
0.1	1.3601696442051\\
0.0888888888888889	1.41360761117434\\
0.08	1.42918189472032\\
0.0727272727272727	1.44048979188526\\
0.0666666666666667	1.46779817772349\\
0.0615384615384615	1.44733187691408\\
0.0571428571428571	1.43491619900104\\
0.0533333333333333	1.45946722120239\\
0.05	1.38217781225949\\
0.0470588235294118	1.44425394416244\\
0.0444444444444444	1.46100042670171\\
0.0421052631578947	1.47131237075794\\
0.04	1.47871418980504\\
0.0298507462686567	1.47066712918519\\
0.0224719101123596	1.42790190081274\\
0.0168067226890756	1.472829512136\\
0.0126182965299685	1.48634509123207\\
0.00947867298578199	1.46811309178901\\
0.00710479573712256	1.48603141822407\\
0.00533333333333333	1.47165959451794\\
0.004	1.47316533642525\\
0.00299850074962519	1.42712945072658\\
0.00224845418774592	0.183270192562159\\
0.00168634064080944	0.0518894920453249\\
0.00126462219411951	0.000124899803876361\\
0.000948541617263457	4.93688992739579e-07\\
0.000711237553342817	3.87766629651125e-09\\
0.00053340445392719	1.30259865900575e-09\\
0.0004	6.59877727072055e-10\\
};

\addplot [color=black, line width=0.5pt, solid, mark=square,mark size=1,mark
options={solid},mark repeat=1]
  table[row sep=crcr]{%
0.2	1.44099439915772\\
0.16	1.48271020626316\\
0.133333333333333	1.42238633970625\\
0.114285714285714	1.32803423407578\\
0.1	1.11612618485511\\
0.0888888888888889	1.14754810983187\\
0.08	0.674716892070011\\
0.0727272727272727	0.53911397796323\\
0.0666666666666667	0.465821547795751\\
0.0615384615384615	0.369791650785775\\
0.0571428571428571	0.317978257762655\\
0.0533333333333333	0.274080327027329\\
0.05	0.255494666524746\\
0.0470588235294118	0.213500813013552\\
0.0444444444444444	0.20294815338449\\
0.0421052631578947	0.178094726776714\\
0.04	0.175344841701929\\
0.0298507462686567	0.115600682795912\\
0.0224719101123596	0.06188049065595\\
0.0168067226890756	0.0266257422624881\\
0.0126182965299685	0.000126829437495462\\
0.00947867298578199	4.29206184408126e-05\\
0.00710479573712256	1.34152560418712e-05\\
0.00533333333333333	4.43693634412123e-06\\
0.004	1.49824053358396e-06\\
0.00299850074962519	4.64697534096611e-07\\
0.00224845418774592	1.47353113284345e-07\\
0.00168634064080944	4.69162638591623e-08\\
0.00126462219411951	1.5046943169078e-08\\
0.000948541617263457	4.94949210591915e-09\\
0.000711237553342817	1.76405927977903e-09\\
0.00053340445392719	8.49634947068195e-10\\
0.0004	6.37758872066709e-10\\
};

\end{axis}
\end{tikzpicture}%
        \input{AccuracyVsCostEx2.tex}\\
        \vspace{0.5cm}
%
\definecolor{mycolor1}{rgb}{0.00000,0.44706,0.74118}%
\definecolor{mycolor2}{rgb}{0.95000,0.32500,0.09800}%
\definecolor{mycolor3}{rgb}{1.00000,0.99000,0.97000}%
\begin{tikzpicture}
    \begin{axis}[%
    width=5.25in, height=0.75in,
    hide axis,
    xmin=0,
    xmax=60,
    ymin=0,
    ymax=0.4,
    legend style={draw=white!15!black,legend cell align=left, legend columns=-1, font=\scriptsize}
    ]

    \addlegendimage{color=colorclassyblue, line width=0.5pt,solid, mark=o, mark
options={solid, colorclassyblue}}
    \addlegendentry{4MaStBMc $\quad$ };
    \addlegendimage{color=red, only marks, mark=asterisk, mark options={solid, red},
mark size = 1}
    \addlegendentry{\ \ \ 4MaStBM\ \ };
    \addlegendimage{color=colorchromeyellow,line width=0.5pt, dashed, mark=*,mark
    size=1,mark options={solid}}
    \addlegendentry{4MaStCC\ \ };
\addlegendimage{color=colorimag, line width=0.5pt, solid}
    \addlegendentry{4MaCC\ \ };
\addlegendimage{color=colorclassyorange, solid, mark=10-pointed star,mark
size=1,mark options={solid}}
    \addlegendentry{6AF10\ \ };
\addlegendimage{color=colorclassyorange, line width=0.75pt, densely dotted,
mark=triangle,mark size=1,mark options={solid}}
    \addlegendentry{6AF5\ \ };
\addlegendimage{color=black, line width=0.5pt, solid, mark=square,mark
size=1,mark options={solid}}
    \addlegendentry{6IKS};

    \end{axis}
\end{tikzpicture}\\
        \vspace{0.5cm}
%
%
\definecolor{mycolor1}{rgb}{0.00000,0.44700,0.74100}%
\definecolor{mycolor2}{rgb}{0.85000,0.32500,0.09800}%
\definecolor{mycolor3}{rgb}{0.92900,0.69400,0.12500}%
\definecolor{mycolor4}{rgb}{0.49400,0.18400,0.55600}%
\definecolor{mycolor5}{rgb}{0.46600,0.67400,0.18800}%
\definecolor{mycolor6}{rgb}{0.30100,0.74500,0.93300}%
\definecolor{mycolor7}{rgb}{0.63500,0.07800,0.18400}%
\begin{tikzpicture}

\begin{axis}[%
width=1.7in, height=2in, scale only axis, xmode=log, xmin=0.0005, xmax=0.2,
xminorticks=true, ymode=log, ymin=1e-9, ymax=10, yminorticks=true, axis
background/.style={fill=white}, xticklabel style =
{xshift=-0.25ex,yshift=-0.1ex}, ylabel={accuracy}, xlabel={time step
\mbox{}}, grid=both]

\addplot [color=colorclassyblue, line width=0.5pt,solid, mark=o, mark
options={solid, colorclassyblue}]
  table[row sep=crcr]{%
0.125	1.82820119003182\\
0.1	0.924836061341598\\
0.0833333333333333	0.280421744378143\\
0.0714285714285714	0.310824705326214\\
0.0625	0.362025805658336\\
0.0555555555555556	0.262277752496006\\
0.05	0.210284316901005\\
0.0454545454545455	0.057652792778427\\
0.0416666666666667	0.0645822852084753\\
0.0384615384615385	0.0262994284182658\\
0.0357142857142857	0.052976383069542\\
0.0333333333333333	0.0126070471460044\\
0.03125	0.00911173910735004\\
0.0294117647058824	0.0183282472105243\\
0.0277777777777778	0.00693466817634361\\
0.0263157894736842	0.0113714317527777\\
0.025	0.00567864111521108\\
0.0186567164179104	0.00164047205849176\\
0.0140449438202247	0.000588987730909764\\
0.0105042016806723	0.000216275860309009\\
0.00788643533123028	7.09705442299766e-05\\
0.00592417061611374	2.35567257388876e-05\\
0.0044404973357016	7.66539784972324e-06\\
0.00333333333333333	2.47467714920772e-06\\
0.0025	7.9344396740619e-07\\
0.00187406296851574	2.54877036190545e-07\\
0.0014052838673412	8.38049098282391e-08\\
0.0010539629005059	2.95464758831499e-08\\
0.000790388871324692	1.23875902148229e-08\\
0.000592838510789661	6.97950451896348e-09\\
0.00044452347083926	5.30249200824753e-09\\
0.000333377783704494	4.77337699044326e-09\\
0.00025	4.69071957905522e-09\\
};

\addplot [color=colorchromeyellow,line width=0.5pt, dashed, mark=*,mark
size=1,mark options={solid},mark repeat=1]
  table[row sep=crcr]{%
0.125	1.82868441224112\\
0.1	0.925290282177921\\
0.0833333333333333	0.280710321631912\\
0.0714285714285714	0.310679070896744\\
0.0625	0.362120914893878\\
0.0555555555555556	0.262337565809272\\
0.05	0.210245784604055\\
0.0454545454545455	0.0576251261342767\\
0.0416666666666667	0.0645630837519469\\
0.0384615384615385	0.0262857099661582\\
0.0357142857142857	0.0529656185232474\\
0.0333333333333333	0.0125993825928535\\
0.03125	0.00910578018666278\\
0.0294117647058824	0.0183232776081513\\
0.0277777777777778	0.00693084470760489\\
0.0263157894736842	0.0113682042250075\\
0.025	0.00567605658053383\\
0.0186567164179104	0.0016397291105165\\
0.0140449438202247	0.000588736973028064\\
0.0105042016806723	0.000216196094390588\\
0.00788643533123028	7.09451308795716e-05\\
0.00592417061611374	2.35486389024996e-05\\
0.0044404973357016	7.66284012085825e-06\\
0.00333333333333333	2.47386481109439e-06\\
0.0025	7.93187096085645e-07\\
0.00187406296851574	2.5479622994826e-07\\
0.0014052838673412	8.37796772614156e-08\\
0.0010539629005059	2.95387590125182e-08\\
0.000790388871324692	1.23853334878098e-08\\
0.000592838510789661	6.97886389084787e-09\\
0.00044452347083926	5.30230698559515e-09\\
0.000333377783704494	4.7733103001279e-09\\
0.00025	4.69069976568501e-09\\
};

\addplot [color=colorclassyorange, solid, mark=10-pointed star,mark
size=1,mark options={solid},mark repeat=1]
  table[row sep=crcr]{%
0.125	1.40514221148514\\
0.1	1.39880004932183\\
0.0833333333333333	1.40201816300814\\
0.0714285714285714	1.39682914944852\\
0.0625	1.15877518485361\\
0.0555555555555556	1.3207695222989\\
0.05	1.35386291233914\\
0.0454545454545455	0.99236893989996\\
0.0416666666666667	0.741069284385327\\
0.0384615384615385	0.695681787275255\\
0.0357142857142857	1.44303217037517\\
0.0333333333333333	1.82526438680505\\
0.03125	1.63040090258374\\
0.0294117647058824	1.52359711803155\\
0.0277777777777778	1.73444088002646\\
0.0263157894736842	1.46719476822369\\
0.025	0.378347268915603\\
0.0186567164179104	0.170352820901986\\
0.0140449438202247	0.0199670294510905\\
0.0105042016806723	0.0519286838794658\\
0.00788643533123028	0.00210857612050879\\
0.00592417061611374	0.00101666738871613\\
0.0044404973357016	0.000323377616140862\\
0.00333333333333333	8.54385209118988e-05\\
0.0025	4.74182936549249e-07\\
0.00187406296851574	2.77483536522498e-07\\
0.0014052838673412	2.24814070864995e-08\\
0.0010539629005059	4.07126299005139e-09\\
0.000790388871324692	4.49145655604136e-09\\
0.000592838510789661	4.54680576119188e-09\\
0.00044452347083926	4.54614273751178e-09\\
0.000333377783704494	4.53513778972628e-09\\
0.00025	4.61449428063052e-09\\
};

\addplot [color=colorclassyorange, line width=0.75pt, densely dotted,
mark=triangle,mark size=1,mark options={solid},mark repeat=1]
  table[row sep=crcr]{%
0.125	1.41421560334881\\
0.1	1.41418823101939\\
0.0833333333333333	1.41456133429578\\
0.0714285714285714	1.4142213083584\\
0.0625	1.41396765238398\\
0.0555555555555556	1.41422428022216\\
0.05	1.41425029949455\\
0.0454545454545455	1.41421356107075\\
0.0416666666666667	1.4142139979383\\
0.0384615384615385	1.41421321043415\\
0.0357142857142857	1.41451076543723\\
0.0333333333333333	1.43565681120671\\
0.03125	1.41549915089193\\
0.0294117647058824	1.41421281644969\\
0.0277777777777778	1.41427856029978\\
0.0263157894736842	1.41768541875128\\
0.025	1.42210610261137\\
0.0186567164179104	1.41421356237448\\
0.0140449438202247	1.41421356200633\\
0.0105042016806723	1.41421356237436\\
0.00788643533123028	1.41421356242206\\
0.00592417061611374	1.41421356237437\\
0.0044404973357016	1.4142137411191\\
0.00333333333333333	0.140054987093893\\
0.0025	3.97122125496013e-05\\
0.00187406296851574	1.2274170277536e-08\\
0.0014052838673412	5.24546696376328e-09\\
0.0010539629005059	4.55307393292485e-09\\
0.000790388871324692	4.56738450833983e-09\\
0.000592838510789661	4.54036823067494e-09\\
0.00044452347083926	4.53892647853882e-09\\
0.000333377783704494	4.53187277521389e-09\\
0.00025	4.61324571460326e-09\\
};

\addplot [color=black, line width=0.5pt, solid, mark=square,mark size=1,mark
options={solid},mark repeat=1]
  table[row sep=crcr]{%
0.125	1.85494302749674\\
0.1	1.35892792944903\\
0.0833333333333333	0.488130494494647\\
0.0714285714285714	0.26497196668039\\
0.0625	0.459677280081579\\
0.0555555555555556	0.380603252191244\\
0.05	0.198323808337455\\
0.0454545454545455	0.0524171094127853\\
0.0416666666666667	0.0447008597341876\\
0.0384615384615385	0.0114680381560209\\
0.0357142857142857	0.0436777561208933\\
0.0333333333333333	0.00214763527462076\\
0.03125	0.00187279773846845\\
0.0294117647058824	0.0100874432412646\\
0.0277777777777778	0.000837311912421154\\
0.0263157894736842	0.00583994129926332\\
0.025	0.000406288999606731\\
0.0186567164179104	6.11789073697043e-05\\
0.0140449438202247	5.98275056956967e-06\\
0.0105042016806723	6.11907841941008e-07\\
0.00788643533123028	1.02235240051469e-07\\
0.00592417061611374	4.37777364948667e-08\\
0.0044404973357016	1.49455166558599e-08\\
0.00333333333333333	7.40374031154423e-09\\
0.0025	5.30796338809436e-09\\
0.00187406296851574	4.76959759369982e-09\\
0.0014052838673412	4.63303159501684e-09\\
0.0010539629005059	4.58366429543514e-09\\
0.000790388871324692	4.57348720372152e-09\\
0.000592838510789661	4.54248469051601e-09\\
0.00044452347083926	4.53968159488868e-09\\
0.000333377783704494	4.53222928927586e-09\\
0.00025	4.61340414162414e-09\\
};

\end{axis}
\end{tikzpicture}%
        \hspace{-0.4cm}
        \input{AccuracyVsCostEx1.tex}\\
        \vspace{0.5cm}
%
\definecolor{mycolor1}{rgb}{0.00000,0.44706,0.74118}%
\definecolor{mycolor2}{rgb}{0.95000,0.32500,0.09800}%
\definecolor{mycolor3}{rgb}{1.00000,0.99000,0.97000}%
\begin{tikzpicture}
    \begin{axis}[%
    width=5.25in, height=0.75in,
    hide axis,
    xmin=0,
    xmax=60,
    ymin=0,
    ymax=0.4,
    legend style={draw=white!15!black,legend cell align=left, legend columns=-1, font=\scriptsize}
    ]

    \addlegendimage{color=colorclassyblue, line width=0.5pt,solid, mark=o, mark
options={solid, colorclassyblue}}
    \addlegendentry{4MaStBMc $\quad$};
    \addlegendimage{color=red, only marks, mark=asterisk, mark options={solid, red},
mark size = 1}
    \addlegendentry{\ \ \ 4MaStBM\ \ };
    \addlegendimage{color=colorchromeyellow,line width=0.5pt, dashed, mark=*,mark
    size=1,mark options={solid}}
    \addlegendentry{4MaStCC\ \ };
\addlegendimage{color=colorimag, line width=0.5pt, solid}
    \addlegendentry{4MaCC\ \ };
\addlegendimage{color=colorclassyorange, solid, mark=10-pointed star,mark
size=1,mark options={solid}}
    \addlegendentry{6AF50\ \ };
\addlegendimage{color=colorclassyorange, line width=0.75pt, densely dotted,
mark=triangle,mark size=1,mark options={solid}}
    \addlegendentry{6AF10\ \ };
\addlegendimage{color=black, line width=0.5pt, solid, mark=square,mark
size=1,mark options={solid}}
    \addlegendentry{6IKS};

    \end{axis}
\end{tikzpicture}
    \caption{[Example 1 (top row), Example 2 (bottom row)] Accuracy {\em vs} time step (left); accuracy {\em vs}
    cost in seconds (right). The prefixes 4 and 6 are used to highlight that these
    are order 4 and order 6 methods, respectively. }
    \label{fig:error}
\end{figure}

In the first example $\ve_1=1$, we use $M_1=150$ spatial grid points while
the highly oscillatory behaviour in the second example, that occurs because
of the small semiclassical parameter $\ve_2=10^{-2}$, requires $M_2=2000$
spatial grid points.

\begin{remark}
Even though IKS behaves like an order four method for $\ve=1$, it is
surprising that it works as well as it does considering that it is designed
for the semiclassical regime $\ve \ll 1$. In particular it outperforms the
sixth-order AF for large time steps.
\end{remark}

{\bf Reference solutions.} The reference solutions for these experiments were
generated by using the sixth-order scheme AF10 with very small time steps. In
both examples we use $M=5000$ grid points and $h=T/10^6$ as time step, which
corresponds to $10^6$ time steps. As usual, the $\CC{L}^2$ distance from the
reference solution is used as a measure of accuracy in these experiments.



\section{Conclusions}
\label{sec:conclusions}
To summarise, in Sections~2 and 3, we have presented a general strategy for
quickly extending any existing fourth-order method for time-independent
potentials to effectively handle the case of laser potentials. The overall
schemes require, at most, one additional exponential, which leads to a very
marginal increase in cost.

Further optimisation steps are possible which allow us to use exactly the
same number of exponentials as the fourth-order methods for time-independent
potentials. In particular, in Sections~2 and 4, we have presented a highly
optimised compact splitting that requires merely five exponentials, requiring
$4$ FFTs (per dimension) per step.

As we can see from the numerical results presented in Figure~\ref{fig:error},
these fourth-order schemes usually exceed the accuracy of the sixth-order
methods of \cite{alvermann2011hocm} for moderate to large time steps, while
being significantly cheaper. Although the specialised sixth-order approach of
\cite{IKS18jcp} has a very high accuracy in the semiclassical regime (as one
might expect), the low costs of our methods makes them highly competitive
with that scheme. The ease of implementing the proposed schemes or extending
existing implementations should make this strategy appealing.

Note that, these methods are not derived under any specific spatial
discretisation choice. The assumption of spectral collocation in
Subsection~\ref{sec:combine} serves only to concretely demonstrate the fact
that certain exponentials can be combined in practice. In principle, it
should be possible to utilise these schemes for strategies such as Hagedorn
wavepackets \cite{gradinaru14nm}, which allow computation over the real line,
or the use of absorbing boundaries, which might be helpful in the case of
non-confining potentials.

\section*{Acknowledgments}
The work of Karolina Kropielnicka in this project was financed by The National Center for
Science, based on grant no. 2016/22/M/ST1/00257.

\bibliographystyle{agsm}
\bibliography{CompactSchemes}

\end{document}